\documentclass[twocolumn]{article}

\usepackage[left=0.6in,right=0.556in,bottom=0.55in,top=0.5in]{geometry}
\usepackage[numbers]{natbib}        
\usepackage[T1]{fontenc}
\usepackage{lmodern}
\usepackage{amsmath,amssymb,amsfonts}
\usepackage{graphicx}
\usepackage{textcomp}
\usepackage{xcolor}
\usepackage{booktabs}
 \usepackage{psfrag}
 \usepackage{caption}
\usepackage{subcaption}
\usepackage{balance}
\usepackage{multirow}
 \usepackage{amsthm}

\usepackage[linesnumbered,ruled,vlined,algo2e]{algorithm2e}

\SetKwInput{KwInput}{Input}                
\SetKwInput{KwOutput}{Output}              

\usepackage{dsfont}
\usepackage{tikz}
\usetikzlibrary{shapes.geometric}
\usepackage{pgfplots}
\pgfplotsset{width=0.45\textwidth}

 \newtheorem{thm}{Theorem}[section]

\usepackage{graphicx}

\newcommand{\T}{^\mathsf{T}} 
\newcommand{\R}{\mathds{R}} 

\usepackage{authblk}
\date{}

\begin{document}

\title{Automatic Scenario Generation for Robust Optimal Control Problems\footnote{This work has received funding from the EPSRC (Engineering and Physical Sciences) under the Active Building Centre project (reference number: EP/V012053/1).  M. Zagorowska also acknowledges funding from the European Research Council (ERC) under the H2020 Advanced Grant no. 787845 (OCAL).}}

%
\author[1]{M. Zagorowska} 
\author[2,5]{P. Falugi} 
\author[3]{E. O'Dwyer}
\author[2,4]{E. C. Kerrigan}

\affil[1]{Automatic Control Laboratory, ETH Zurich, Switzerland (e-mail: mzagorowska@ethz.ch)}
\affil[2]{Department of Electrical and Electronic Engineering, Imperial College London, UK (e-mail: p.falugi@imperial.ac.uk, e.kerrigan@imperial.ac.uk).}
\affil[3]{Department of Chemical Engineering, Imperial College London, UK  (e-mail: e.odwyer@imperial.ac.uk)}
\affil[4]{Department of Aeronautics, Imperial College London, UK}
\affil[5]{Department of Engineering \& Construction, University of East London, UK.}
\maketitle
\begin{abstract}                
Existing methods for nonlinear robust control often use scenario-based approaches to formulate the control problem as nonlinear optimization problems. Increasing the number of scenarios improves robustness, while increasing the size of the optimization problems. Mitigating the size of the problem by reducing the number of scenarios requires knowledge about how the uncertainty affects the system. This paper draws from local reduction methods used in semi-infinite optimization to solve robust optimal control problems with parametric uncertainty. We show that nonlinear robust optimal control problems are equivalent to semi-infinite optimization problems and can be solved by local reduction. By iteratively adding interim globally worst-case scenarios to the problem, methods based on local reduction provide a way to manage the total number of scenarios. In particular, we show that local reduction methods find worst case scenarios that are not on the boundary of the uncertainty set. The proposed approach is illustrated with a case study with both parametric and additive time-varying uncertainty. The number of scenarios obtained from local reduction is 101, smaller than in the case when all $2^{14+3\times192}$ boundary scenarios are considered. A validation with randomly drawn scenarios shows that our proposed approach reduces the number of scenarios and ensures robustness even if local solvers are used.
\end{abstract}



\section{Introduction}
Robust nonlinear optimal control problems are often solved using a scenario-based approach, where each \emph{scenario} corresponds to a separate realization of uncertainty. This paper draws from approaches used in semi-infinite optimization to solve robust optimal control problems in an efficient way.

To ensure that the optimization problems resulting from scenario-based approaches to robust control are tractable, the number of scenarios must be limited \citep{scenario_Calafiore2006}. Usually, the choice of scenarios is done from experience \citep{scenario_Grammatico2015,Robust_Puschke2018} and requires knowledge about both the controlled system and the uncertainty to ensure that the chosen scenarios guarantee robustness. A recent review of scenario-based methods was done by \cite{scenario_Campi2021} who indicated that scenario selection is highly affected by the knowledge about the uncertainty distribution. In practice, to tackle problems with limited knowledge about the uncertainty, it is often assumed that the worst-case scenarios lie on the boundary of the uncertainty set \citep{Cutting_Mutapcic2009,Handling_Lucia2014,Monotonicity_Vuffray2015}. As indicated by \cite{Sensitivity_Thombre2021}, the worst-case scenario in nonlinear systems may lie in the interior of the uncertainty range. In this paper, we  present a method for choosing potential worst-case scenarios that is derived from semi-infinite optimization and is independent from the uncertainty distribution. 

An in-depth review of semi-infinite optimization methods has been done by \cite{Semi_Hettich1993,Semi_Hettich2009,review_Hettich1983} and recently by  \cite{adaptive_Seidel2020,Recent_Djelassi2021}.  Semi-infinite optimization methods have been used for optimal control by \cite{Semi_Hauser2018} to find optimal trajectories for robotic arms. However, they considered only exogenous uncertainty due to obstacles that did not affect the dynamics of the controlled systems. \cite{Robust_Puschke2018} used semi-infinite optimization methods to solve an optimal control problem with only parametric uncertainty. Time-varying uncertainty was considered by \cite{Sensitivity_Thombre2021} who used local reduction to find the interim worst-case scenarios. However, they assumed that at every time step the number of possible scenarios was finite.

The main contribution of the current work is a formulation of robust nonlinear optimal control problems as semi-infinite optimization problems. We then numerically demonstrate that local reduction methods enable more robust handling of significant parametric and time-varying uncertainty than existing approaches.

The rest of the paper is structured  as follows. Section~\ref{sec:ProblemFormulation} introduces robust optimal control problems. Section~\ref{sec:LocalReduction} presents the new method for solving robust optimal control problems. The numerical results are shown in Section~\ref{sec:NumericalResults}. The paper ends with conclusions in Section~\ref{sec:Conclusions}.

\section{Problem formulation}
\label{sec:ProblemFormulation}

\subsection{Semi-infinite optimization problem}
\label{sec:semiinfProblem}
\begin{subequations} \label{eqn:SemiInf}
\begin{align}
\mathcal{Q}:\quad &\min_{\theta\in\mathcal{A}} \quad Q(\theta)
    \label{eqn:SemiinfCost}\\
&\text{subject to } R(\theta,\rho)\leq 0 \label{eq:Semiinfmapping} \text{ for all }\rho\in\mathcal{B}
\end{align}
\end{subequations}
where $\mathcal{A}\subset \R^{n_{\theta}}$ and $\mathcal{B}\subset\R^{n_{\rho}}$ are nonempty and compact sets, and $Q$ and $R$ are continuous functions of their respective arguments \citep{Infinitely_Blankenship1976}. The problem \eqref{eqn:SemiInf} has a finite number of variables $\theta$ but includes an infinite number of constraints if $\mathcal{B}$ has an infinite number of points. In particular, $\mathcal{B}$ may be uncountable. 

One approach to remove the infinite number of constraints consists in rewriting the constraint \eqref{eq:Semiinfmapping} as:
\begin{equation}
    S(\theta):=\max_{\rho\in\mathcal{B}} R(\theta,\rho)\leq 0. \label{eq:SemiinfmappingMax}
\end{equation}
The challenge in solving the equivalent problem with constraint \eqref{eq:SemiinfmappingMax} is in non-differentiability of the function $S(\cdot)$. The local reduction method proposed by \cite{Infinitely_Blankenship1976} allows overcoming the non-differentiability of $S(\cdot)$ by sequentially solving \eqref{eqn:SemiInf} with finite subsets of constraints taken from $\mathcal{B}$. In this paper, we show that optimal control problems can be formulated as semi-infinite optimization problems and solved using the method proposed by \cite{Infinitely_Blankenship1976}.

\subsection{Dynamic system with uncertainty}
The system to be controlled is described by a nonlinear difference equation with time-varying uncertainty $w_k\in\mathbb{W}\subset \R^{n_w}$ and constant uncertainty $d\in\mathbb{D}\subset\R^{n_d}$:
\begin{equation}
x_{k+1}=f_k(x_k,u_k,w_k,d)
\label{eq:InitialDynamics}
\end{equation}
where $f_k$ is continuously differentiable.
The state $x_0$ at time zero is w.l.o.g.\ assumed to be equal to a given~$\hat{x}$. 

The control trajectory $\mathbf{u}:=(u_0,\ldots,u_{N-1})$ is generated by a causal dynamic feedback policy  
\[
u_k:=\pi_k(x_0,\ldots,x_k;q_0,\ldots,q_k,r)
\] that is parameterised by $\mathbf{q}:=(q_0,q_1,\ldots,q_{N-1})\in\mathbb{R}^{n_q}$ and $r\in\R^{n_r}$. The state trajectory $\mathbf{x}:=(x_0,\ldots,x_N)$. The time-varying uncertainty $w_k$ at time~$k$ and the constant uncertainty $d$ affect the dynamics in both an additive and non-additive way, and take on values from compact and uncountable (infinite cardinality) sets. Uncertainty in the measured value of $x_k$ can be modelled by a suitably-defined choice of $f_k$, $\pi_k$ and $w_k$.

A trajectory $(\mathbf{x},\mathbf{u})$ satisfying the dynamics \eqref{eq:InitialDynamics} and control policy for a given  parameterization $(\mathbf{q},r)$ and realisation of uncertainty $(\mathbf{w}$, $d$), where the trajectory $\mathbf{w}:=(w_0,\ldots,w_{N-1}) \in \mathbb{W}^N:=\mathbb{W}\times\cdots\times\mathbb{W}$,  is defined as:
\begin{equation}
\label{eq:Trajectory}
\begin{aligned} \mathbf{z}(\mathbf{q},r,\mathbf{w},d):=  \Big\lbrace (\mathbf{x},\mathbf{u}) &\mid x_0=\hat{x}\\
x_{k+1}&=f_k(x_k,u_k,w_k,d)\\
u_k&=\pi_k(x_0,\ldots,x_k;q_0,\ldots,q_k,r)\\
k&=0,1,\ldots,N-1 \Big\rbrace. \end{aligned} 
\end{equation}

\subsection{Robust optimal control problem}
\label{sec:RobustOptimalControl}

\subsubsection{Objective function and constraints}
\label{sec:ObjCstr}
The cost function for the optimal control problem over a horizon of length $N$ is:
\begin{equation}
J_N(\mathbf{x},\mathbf{u},\mathbf{w}, d):=J_f(x_N,w_N,d)+\sum\limits_{k=0}^{N-1}\ell_k(x_k,u_k,w_k,d).
\label{eq:Objective}
\end{equation}
Both the terminal cost function $J_f(\cdot,\cdot,\cdot)$ and stage cost  $\ell_k(\cdot,\cdot,\cdot,\cdot)$ are continuously differentiable and depend on the uncertainty $w$ and $d$. The objective of the optimal control problem is to find a feedback policy $\pi$ for system~\eqref{eq:InitialDynamics} such that the worst-case cost in \eqref{eq:Objective} is minimized and the constraints
\begin{equation}
g_k(x_k,u_k,w_k,d)\leq 0
\label{eq:InitCstr}
\end{equation}
are satisfied for all time instants $k=0,\ldots,N-1$, all states~$\mathbf{x}$, control~$\mathbf{u}$, uncertainty~$\mathbf{w}$~and~$d$. The vector function of $n_{g}$ components, $g_k(\cdot,\cdot,\cdot,\cdot)$, is continuously differentiable and depends on uncertainty $\mathbf{w}$ and $d$. Note that a constraint on~ $x_N$ can be included by incorporating $f_{N-1}$ in a suitable definition of $g_{N-1}$.

\subsubsection{Semi-infinite formulation}
\label{sec:SemiInfForm}

Given a set of uncertainties 
$
\mathbb{H} \subseteq \mathbb{W}^N\times\mathbb{D},
$
the problem in this work is stated as:
\begin{subequations}
\label{eq:Plifted}
\begin{align}
 \mathcal{P}_N(\mathbb{H}): \min_{{\substack{\mathbf{q},r\\ \mathbf{x}^i,\mathbf{u}^i,i\in\mathbb{J}}}}& \max_{i\in\mathbb{J}}\  J_N(\mathbf{x}^i,\mathbf{u}^i,\mathbf{w}^i, d^i) \label{eq:CostLifted}
   \\
 \text{s.t. } 
g_k(x_k^i,u_k^i,w_k^i,d^i)\leq 0,\ &\forall i\in\mathbb{J}, k=0,\ldots,N-1\\
(\mathbf{x}^i,\mathbf{u}^i) = \mathbf{z}(\mathbf{q},r,\mathbf{w}^i,d^i),\ &\forall i\in\mathbb{J} \label{eq:TrajectoryAll}
\end{align}
\end{subequations}
where $\mathbb{J} := \lbrace 1,\ldots,\operatorname{card} \mathbb{H}\rbrace$ and $(\mathbf{x}^i,\mathbf{u}^i)$ is the state and input trajectory associated with the $i^\text{th}$ disturbance realisation $(\mathbf{w}^i,d^i)$ such that $\mathbb{H} = \bigcup_{i\in\mathbb{J}}\{(\mathbf{w}^i,d^i)\}.$

If $\mathbf{z}(\cdot)$ in \eqref{eq:TrajectoryAll} is linear jointly in all arguments, the problem \eqref{eq:Plifted} can often be solved using scenario-based methods for robust control from \cite{scenario_Calafiore2006,Min_Scokaert1998}, provided additional convexity assumptions are satisfied by the uncertainty set $\mathbb{W}$. In this work, the dynamics from \eqref{eq:InitialDynamics} are nonlinear and $\mathbb{W}$ is only non-empty and compact.

\begin{thm}
\label{thm:MainResult}
The robust optimal control problem \eqref{eq:Plifted} is equivalent to the semi-infinite optimization problem \eqref{eqn:SemiInf} with $\theta:=(\mathbf{q},r,\gamma)$, where $\gamma$ is an additional scalar parameter characterizing the cost upper-bound, $\rho:=(\mathbf{w},d)$ and the sets $\mathcal{A}:=\R^{n_q}\times\R^{n_r}\times\R$, $\mathcal{B}:=\mathbb{H}$.
\end{thm}
\begin{proof}
In contrast to \eqref{eqn:SemiinfCost}, the objective function in \eqref{eq:CostLifted} contains uncertainty. Introducing $\gamma\in\R$, we rewrite \eqref{eq:Plifted} as:
\begin{subequations}
\label{eq:Pgamma}
\begin{align}
 \mathcal{P}_N(\mathbb{H}):\quad \min_{{\substack{\gamma,\mathbf{q},r\\ \mathbf{x}^i,\mathbf{u}^i, i\in\mathbb{J}}}}& \quad \gamma 
   \\
 \text{s.t. } 
g_k(x_k^i,u_k^i,w_k^i,d^i)\leq 0,\ &\forall i\in\mathbb{J},k=0,\ldots,N-1 \label{eq:PaugIneq}\\
(\mathbf{x}^i,\mathbf{u}^i) = \mathbf{z}(\mathbf{q},r,\mathbf{w}^i,d^i),\ &\forall i\in\mathbb{J}\\
J_N(\mathbf{x}^i,\mathbf{u}^i,\mathbf{w}^i, d^i)\leq \gamma,\ &\forall i\in\mathbb{J}\label{eq:CstrObjRef}
\end{align}
\end{subequations}
The problem \eqref{eq:Pgamma} has uncertainty exclusively in the constraints. If $\operatorname{card} \mathbb{H}$ is finite, then the problem \eqref{eq:Pgamma} is convenient to solve numerically using tailored efficient finite-dimensional optimization methods that exploit the sparsity in the relevant Jacobians and Hessians. However, infinite cardinality of $ \mathbb{H}$ yields an infinite number of both constraints and variables, which means that the problem \eqref{eq:Pgamma} needs to be further reformulated to become \eqref{eqn:SemiInf}. Noticing that the constraint \eqref{eq:PaugIneq} is equivalent to
\begin{equation}
    \max_k g_k(x_k^i,u_k^i,w_k^i,d^i)\leq 0,\ \forall i\in\mathbb{J},
\end{equation}
we introduce
\begin{multline}
G(\mathbf{x}^i,\mathbf{u}^i,\mathbf{w}^i,d^i,\gamma):=\max\lbrace \max_{h,k}\; e_h^T g_k(x^i_k,u^i_k,w^i_k,d^i),\\ J_N(\mathbf{x}^i,\mathbf{u}^i,\mathbf{w}^i,d^i)-\gamma \rbrace
\label{eq:AllConstraints}    
\end{multline}
In \eqref{eq:AllConstraints}, $e_h$ is the $h^\text{th}$ column of an identity matrix $\mathbb{I}_{n_g} $. Using~\eqref{eq:Trajectory} and \eqref{eq:AllConstraints}, we can write \eqref{eq:Pgamma} as:
\begin{subequations}
\label{eq:Prevised}
\begin{align}
&\mathcal{P}_N(\mathbb{H}):\quad  \min_{\mathbf{q},r,\gamma} \quad \gamma \label{eq:RevisedCost}
   \\
 &\quad \text{s.t. } 
G(\mathbf{z}(\mathbf{q},r,\mathbf{w},d),\mathbf{w},d,\gamma)\leq 0,\ \forall (\mathbf{w},d)\in\mathbb{H},
\label{eq:PreConstraints}
\end{align}
\end{subequations}
The problem \eqref{eq:Prevised} is equivalent to
\begin{subequations}
\label{eq:PrevisedMax}
\begin{align}
&\mathcal{P}_N(\mathbb{H}):\quad  \min_{\mathbf{q},r,\gamma} \quad \gamma
   \\
 &\quad \text{s.t. } 
\max_{(\mathbf{w},d)\in\mathbb{H}}G(\mathbf{z}(\mathbf{q},r,\mathbf{w},d),\mathbf{w},d,\gamma)\leq 0.
\label{eq:PreConstraintsMax}
\end{align}
\end{subequations}

Taking $\theta:=(\mathbf{q},r,\gamma)$ and $\rho:=(\mathbf{w},d)$ in \eqref{eq:Prevised} (similarly in \eqref{eq:PrevisedMax}) we obtain the form of \eqref{eqn:SemiInf} (similarly \eqref{eq:SemiinfmappingMax}).\hspace*{\fill}\qed
\end{proof}

We also notice that \eqref{eq:PreConstraintsMax} is equivalent to:
\begin{equation}
    G_{\max}(\mathbf{q},r,\gamma,\mathbb{H}):=\max_{
{\substack{(\mathbf{w},d)\in\mathbb{H}\\ (\mathbf{x},\mathbf{u})=\mathbf{z}(\mathbf{q},r,\mathbf{w},d)}}
} G(\mathbf{x},\mathbf{u},\mathbf{w},d,\gamma)
    \label{eq:Gmax}
\end{equation}

Theorem~\ref{thm:MainResult} allows solving problem~\eqref{eq:Plifted} as a semi-infinite optimization problem of the form~\eqref{eqn:SemiInf} using the local reduction from \cite{Infinitely_Blankenship1976}.

\section{Local reduction for optimal control}
\label{sec:LocalReduction}

The local reduction method from \cite{Infinitely_Blankenship1976} consists in iteratively solving finite-dimensional optimization problems. We use the local reduction methods for the problem  \eqref{eq:Prevised} by iteratively solving optimal control problems parametrised by scenarios. A scenario is a realisation of the uncertainty $(\mathbf{w}^*,d^*)\in\mathbb{W}^N\times\mathbb{D}$. 

\subsection{Minimization step}
Algorithm~\ref{alg:LocalReductionGeneral} is the local reduction algorithm for robust optimal control. The  algorithm in iteration $j$ solves an optimal control problem of the form \eqref{eq:Pgamma} or \eqref{eq:Prevised} assuming that the number of scenarios $\operatorname{card} \mathbb{H}_j$ at step~$j$ is finite. The algorithm needs an initial guess for the parameters of the controller. For instance, the initial guess can be obtained by solving \eqref{eq:Prevised} for one scenario, $\operatorname{card}\mathbb{H}_1=1$.

In the first step of Algorithm~\ref{alg:LocalReductionGeneral} (line 3), the algorithm checks whether worst-case scenarios exist that would lead to a violation of constraints~\eqref{eq:AllConstraints}. If no constraints are violated (line 4), the current parameters give a robust solution to the current set of scenarios $\mathbb{H}_j$. If there exists at least one violated constraint, then a scenario corresponding to the maximum constraint violation is added to the scenario set $\mathbb{H}_{j+1}$ in the next iteration (line 7). The new set $\mathbb{H}_{j+1}$  is then used to find a new set of control parameters (line~9). The algorithm ends if no new scenarios are added, i.e.\ $\operatorname{card}\mathbb{H}_{j}= \operatorname{card}\mathbb{H}_{j-1}$. 

In this work, any scenario corresponding to the maximum constraint violation can be added to the set of scenarios. However, it has been shown that computational performance may be improved if multiple scenarios are added \citep{global_Tsoukalas2008}.

\setlength{\algomargin}{1.1em}
\begin{algorithm2e}[t]
\SetAlgoLined
  \KwInput{Initial guess for $\mathbf{q}$, $r$, $\gamma$ and $\mathbb{H}_1\neq\emptyset$}
  \KwOutput{Optimal  $\mathbf{q}^*$, $r^*$, $\gamma^*$, set of scenarios $\mathbb{H}^*$}
  Set $\mathbf{q}^1\leftarrow\mathbf{q}$, $r^1\leftarrow r$, $\gamma^1\leftarrow\gamma$, $j\leftarrow 1$
  
  \Repeat{$\operatorname{card}\mathbb{H}_{j}= \operatorname{card}\mathbb{H}_{j-1}$}{
Compute $G_{\max}(\mathbf{q}^j,r^j,\gamma^j,\mathbb{W}^N\times\mathbb{D})$
and a maximizer $(\mathbf{x}^j,\mathbf{u}^j,\mathbf{w}^j,d^j)$ by solving~\eqref{eq:Gmax} with $\mathbb{H}=\mathbb{W}^N\times \mathbb{D}$.

\eIf{$G_{\max}(\mathbf{q}^j,r^j,\gamma^j,\mathbb{W}^N\times\mathbb{D})\leq 0$}{
\begin{equation*}
  \mathbb{H}_{j+1}\leftarrow\mathbb{H}_{j}
  \end{equation*}  }{
  Add new scenario 
  \begin{equation*}
  \mathbb{H}_{j+1}\leftarrow\mathbb{H}_{j}\cup \{(\mathbf{w}^j,d^j)\}
  \end{equation*}
      
  Find a $(\mathbf{q}^{j+1},r^{j+1},\gamma^{j+1})$ that solves $\mathcal{P}_N(\mathbb{H}_{j+1})$ using~\eqref{eq:Pgamma} or \eqref{eq:Prevised}.
}
  Set $j\leftarrow j+1$

      Set $(\mathbf{q}^*,r^*,\gamma^*)\leftarrow(\mathbf{q}^j,r^j,\gamma^j)$  and $\mathbb{H}^*\leftarrow \mathbb{H}_j$.
 }
 \caption{Exact local reduction method\label{alg:LocalReductionGeneral}}
\end{algorithm2e}

The convergence of local reduction method in the case of the form \eqref{eqn:SemiInf} was first shown by \cite{Infinitely_Blankenship1976}, generalised by \cite{Semi_Reemtsen1998}, and recently by \cite{Global_Mitsos2011}. They required that the sets $\mathcal{A}$ and $\mathcal{B}$ in \eqref{eqn:SemiInf} are non-empty and compact, and that the functions $Q$ and $R$ are continuous with respect to all their arguments. Then they showed that the sequence of solutions obtained for a sequence of finite and countable subsets of $\mathcal{B}$ converges to the solution of \eqref{eqn:SemiInf}, provided that the minimization and the maximization steps are solved to global optimality. We show in Theorem \ref{thm:Convergence} when the Algorithm \ref{alg:LocalReductionGeneral} solves the problem \eqref{eq:Plifted}.  A discussion on convergence rate of local reduction methods was done by \cite{adaptive_Seidel2020}. 

\begin{thm}
\label{thm:Convergence}
The solution $(\mathbf{q}^*,r^*,\gamma^*)$ obtained from Algorithm \ref{alg:LocalReductionGeneral} for a non-empty and compact set $\mathbb{W}^N\times\mathbb{D}$ converges to the solution of \eqref{eq:Plifted} if the set $\mathbb{F}$ such that $(\mathbf{q},r,\gamma)\in\mathbb{F}\subset \R^{n_q}\times\R^{n_r}\times\R$ is non-empty and compact.
\end{thm}
\begin{proof}
From Theorem \ref{thm:MainResult}, we have $\theta:=(\mathbf{q},r,\gamma)$ and $\rho:=(\mathbf{w},d)$, $\mathcal{A}:=\mathbb{F}$, $\mathcal{B}:=\mathbb{H}$. In \eqref{eq:Prevised}, we take $Q(\theta):=\gamma$ which is linear and thus continuous. Then we have $R(\theta,\rho):=G(\mathbf{z}(\mathbf{q},d,\mathbf{w},d),\mathbf{w},d)$ which is continuous because both $\max_{h,k} e_h^{T}g_k(\cdot,\cdot,\cdot,\cdot)$ and $J_N(\cdot,\cdot,\cdot,\cdot)$ are continuous. The proof follows from Lemma 2.2 by \cite{Global_Mitsos2011}.\hspace*{\fill}\qed
\end{proof}

\subsection{Maximization step}
The maximization step consists in solving \eqref{eq:Gmax} with $\mathbb{H}=\mathbb{W}^N\times\mathbb{D}$. Solving \eqref{eq:Gmax} is equivalent to solving $n_{g}\cdot (N-1)+1$ optimization problems, where $n_g$ denotes the number of elements in the vector function $g(\cdot)$ from constraints in \eqref{eq:InitCstr}. The algorithm is presented in Algorithm~\ref{alg:LocalReductionMaximisation}. Without loss of generality, we assume that the first constraint to include in the maximization problem corresponds to the reformulated objective function \eqref{eq:Objective}. A scenario that corresponds to maximal value of this constraint is added to an auxiliary set $\mathbb{K}$. The remaining $n_g\cdot (N-1)$ constraints are included as objectives in the respective maximization problems (line four to eight in Algorithm~\ref{alg:LocalReductionMaximisation}). Note that the problem corresponding to the objective (line 2) and all the problems corresponding to the constraints (line four to eight) can be solved in parallel.

All maximization problems are subject to the same equality constraints capturing the dynamics. This formulation allows us to treat the maximization problems as optimal control problems and preserve the sparsity of the relevant Jacobians and Hessians. We solve the maximization problems as optimal control problems where $\mathbf{q}$, $r$, and $\gamma$ are known parameters whereas $\mathbf{w}$ and $d$ are treated as unknown inputs. Thus, the maximization problems can be solved using off-the-shelf optimal control solvers.

Solving~\eqref{eq:Gmax} with $\mathbb{H}=\mathbb{W}^N\times\mathbb{D}$ corresponds to lines four to eight in Algorithm~\ref{alg:LocalReductionMaximisation} and can be done by solving a number of finite-dimensional optimization problems in parallel as indicated by \cite{parallel_Zakovic2003}.

\begin{algorithm2e}[t]
\SetAlgoLined
  \KwInput{Current values of $\mathbf{q}^j$, $r^j$, $\gamma^j$}
  \KwOutput{Worst case scenario $(\mathbf{w}^{j},d^{j})$ in iteration $j$}
  
        Find any $\mathbf{x}^*,\mathbf{u}^*,\mathbf{w}^*,d^*$   that solves:
        \begin{subequations}
        \begin{align}
        \max_{\mathbf{x},\mathbf{u},\mathbf{w},d} & \quad J_N(\mathbf{x},\mathbf{u},\mathbf{w},d)-\gamma^j\nonumber   \\
         \quad \text{s.t. } &
        (\mathbf{x},\mathbf{u}) = \mathbf{z}(\mathbf{q}^j,r^j,\mathbf{w},d)\nonumber\\
        & (\mathbf{w},d) \in \mathbb{W}^N\times\mathbb{D}\nonumber
        \end{align}
        \end{subequations}
        
        Set $\mathbb{K}\leftarrow{(\mathbf{w}^*,d^*,J_N(\mathbf{x}^*,\mathbf{u}^*,\mathbf{w}^*,d^*)-\gamma^j)}$

  \For{$h=1,\ldots,n_g$}{
    \For{$k=1,\ldots,N-1$}{
          Find any $\mathbf{x}^*,\mathbf{u}^*,\mathbf{w}^*,d^*$  that solves:
        \begin{subequations}
\begin{align}
\max_{\mathbf{x},\mathbf{u},\mathbf{w},d} & \quad e_h g_k(x_k,u_k,w_k,d)\nonumber   \\
 \quad \text{s.t. } &
(\mathbf{x},\mathbf{u}) = \mathbf{z}(\mathbf{q}^j,r^j,\mathbf{w},d)\nonumber\\
        & (\mathbf{w},d) \in \mathbb{W}^N\times\mathbb{D}\nonumber
\end{align}
\end{subequations}
        Set $\mathbb{K}\leftarrow\mathbb{K}\cup{(\mathbf{w}^*,d^*,e_h g_k(x_k^*,u_k^*,w_k^*,d^*))}$
        
        }
    }
    Set $v^*\leftarrow \max \lbrace v_3 \mid (v_1,v_2,v_3)\in\mathbb{K}\rbrace$;
    
    Choose any $(\mathbf{w}^{j},d^{j})\in \lbrace (v_1,v_2) \mid (v_1,v_2,v^*)\in \mathbb{K} \rbrace$
     \caption{Maximization - line 3 in Algorithm~\ref{alg:LocalReductionGeneral} \label{alg:LocalReductionMaximisation}}
\end{algorithm2e}

\section{Numerical results}
\label{sec:NumericalResults}
We show that the local reduction method described in Section~\ref{sec:LocalReduction} finds scenarios from inside the uncertainty sets and provides robust solutions to optimal control problems with uncertainty in a case study of temperature control in a residential building, adapted from \cite{System_Lian2021}. The example was implemented in Julia 1.6 \citep{Julia_Bezanson2017} using JuMP 0.21.4 \citep{JuMP_Dunning2017}. The problems were then solved with Ipopt version 3.12.10 \citep{Exact_Thierry2020}. All tests were performed on a laptop with an Intel\textsuperscript{\textregistered} Core\textsuperscript{\texttrademark} i7-7500U with 16\,GB of RAM.

\subsection{Example} The example shows a single zone building affected by time-varying internal heat gain, solar radiation, and external temperature. The objective is to keep the internal temperature $x_k^{\text{temp}}$ within time-varying bounds:
\begin{equation}
T_{\min}\leq x_k^{\text{temp}}\leq T_{\max}.
\label{eq:TempBounds}
\end{equation}
During the day, the indoor temperature must be kept above $23\,^\circ$C and during the night it can drop down to $17\,^\circ$C.  The maximal temperature is the same during the day and night, $T_{\max}=26\,^{\circ}$C.
The dynamics are discrete and linear:
\begin{equation}
x_{k+1}=A(d)x_k+B(d)u_k^{\text{sat}}+Ww_k.
\label{eq:BuildingSimple}
\end{equation}
The states $x$ describe the indoor temperature~$x^{\text{temp}}$, wall temperature~$x^{\text{wall}}$, and the corridor temperature~$x^{\text{corr}}$. The initial condition was chosen as $$x_0={\begin{bmatrix} 25\,^{\circ}\text{C}& 24\,^{\circ}\text{C} &24\,^{\circ}\text{C} \end{bmatrix}^{\T}}.$$ The control~$u$ represents the amount of heating and cooling delivered to the building. The nominal matrices are taken from \cite{System_Lian2021}. 

There are 14 uncertain parameters affecting the matrices $A$, $B$, and the initial condition for the wall and corridor temperatures, $d\in\R^{14}$, $d=[\lambda_{\text{wall}},\lambda_{\text{corr}},\delta_{i,j},\eta_j]_{i,j=1,\ldots,3}$. We assume that the wall temperature and the corridor temperature can only be measured approximately. We have $x_0^i=24+\lambda_i$, $i=\text{wall}, \text{corr}$, where $\lambda_i\in[-0.5,0.5]$. We also assume that the matrices $A=[a_{i,j}]$ and $B=[b_j]$, $i,j=1,2,3$ are affected by uncertainty $a_{i,j}\cdot\delta_{i,j}$ and $b_j\cdot\eta_j$ where $\delta_{i,j},\eta_j$ are uncertain parameters, $\delta_{i,j},\eta_j\in[0.96,1.03]$.
 
The minimal control effort is ensured by the objective function $$J=\frac{1}{N}\sum\limits_{k=1}^{N}u_k^2.$$ It is assumed that the day starts at 6.00\,am and lasts 12~hours.  The optimal control problem is solved over a period of 48~hours starting at 6.00\,am the first day, with $N=192$. The three uncertain parameters, internal gain, solar radiation, and external temperature, vary with time within the limits provided by \cite{System_Lian2021}. 

The control variables are parameterised as:
\begin{equation}
u_k=Kx_k^{\text{temp}}+q_k
\label{eq:ControlBuilding}
\end{equation}
where $r:=\operatorname{vec}(K)$ and $q_k$ are decision variables. We include saturation of the control inputs:
\begin{equation}
u^{\text{sat}}_k = \begin{cases}
-500\text{\,W}&\text{for } u_k<-500\text{\,W}\\
u_k&\text{for } -500\text{\,W}\leq u_k\leq1200\text{\,W}\\
1200\text{\,W}&\text{for } u_k>1200\text{\,W}
\end{cases}
\end{equation}
The saturation was approximated by a smooth function:
\begin{equation}
u^{\text{sat}}_k=\frac{\beta_0}{\beta_1+\exp(\beta_2u_k)}+\beta_3
\label{eq:SmoothSaturation}
\end{equation}
where $\beta_i$ are constants. Here $\beta_0=-5030$, $\beta_1=2.937$, $\beta_2=0.003$, $\beta_3=1207$.

The results obtained from local reduction are then compared with three scenario-based approaches from the literature \citep{Sensitivity_Thombre2021}: nominal, with a controller obtained assuming there is no uncertainty (``Nominal''), randomised, with a controller obtained for a number of scenarios chosen from a uniform distribution (``Random''), extreme, with a controller obtained for three cases: nominal, lower bound, and upper bound for all uncertainties (``Nominal+two extreme'').

\subsection{Overall performance}
\begin{figure}
     \centering
         \includegraphics[width=0.4\textwidth]{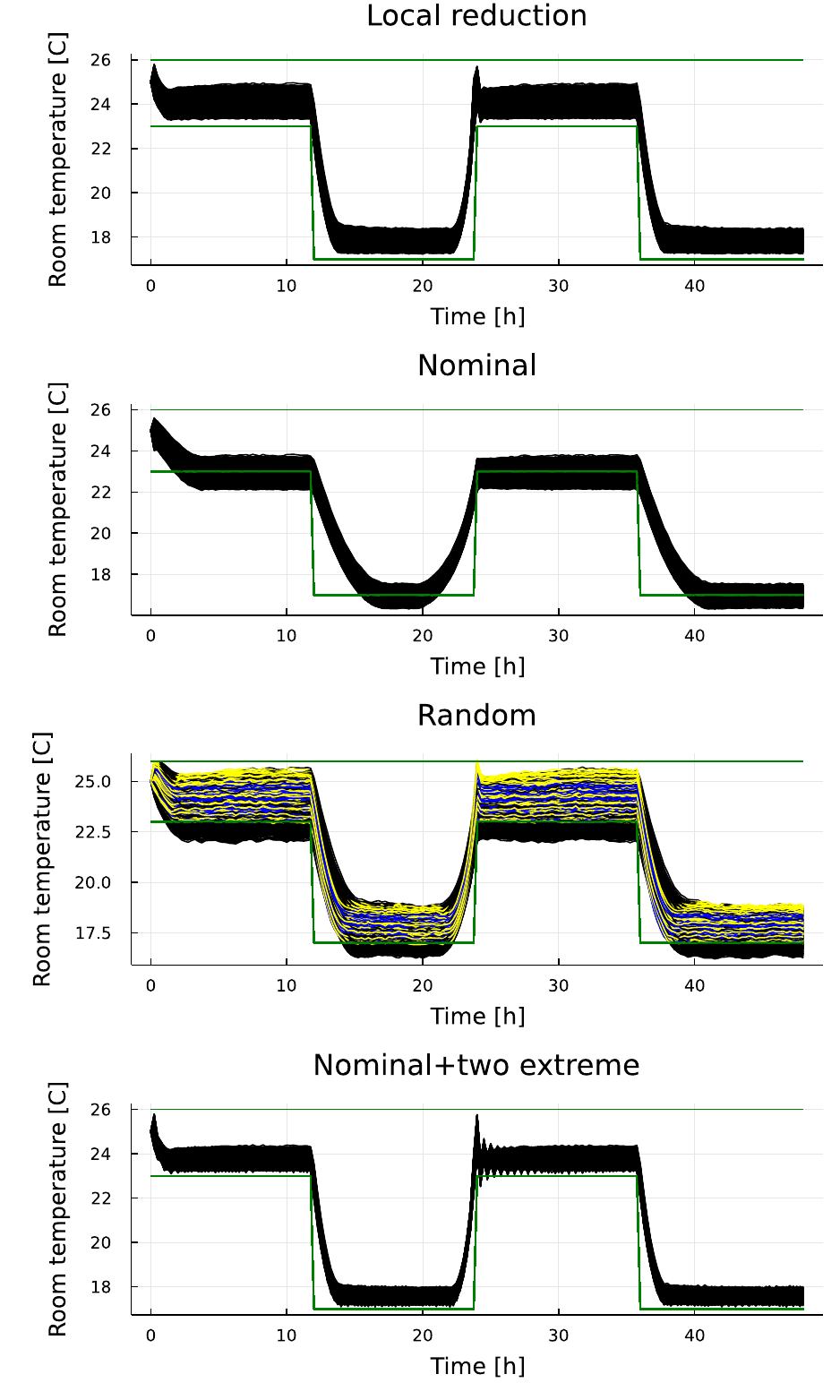}
        \caption{Comparison with other scenario-based approaches}
        \label{fig:Comparison}
\end{figure}

The local reduction method reduced the number of scenarios to 101. The resulting controller obtained for the interim worst-case scenarios was then validated for 500 random realisations from a uniform distribution of uncertainty. The validation of the controller is shown in top plot in Fig.~\ref{fig:Comparison}. The black curves stay within the green bounds corresponding to constraints \eqref{eq:TempBounds}. The results suggest that local reduction was able to find a robust solution despite using a local solver for maximizations.

\begin{figure}
     \centering
         \includegraphics[width=0.33\textwidth]{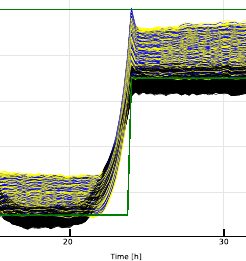}
        \caption{Zoomed in part showing how random scenarios violate constraints }
        \label{fig:ComparisonClipped}
\end{figure}

The results also indicate that the local reduction method handles time-varying uncertainty without specifying the scenarios over a shorter time horizon. Algorithm 2  treats a scenario of time-varying uncertainty as a single realization of the uncertainty over the whole horizon, removing the necessity of defining all possible uncertainty realizations on short horizons. 

\subsection{Comparison with other approaches}

Validating the nominal controller with 500 random scenarios shows that the approach based on nominal values leads to violation of constraints (plot `Nominal' in Fig.~\ref{fig:Comparison}). The second set of controllers we used was derived using three random cases: five scenarios, 100 scenarios, and 250 scenarios. The results are shown in the plot `Random' in Fig.~\ref{fig:Comparison}, with black corresponding to the controller obtained from five scenarios, yellow to the controller with 100 scenarios, and blue to the controller with 250 scenarios. As shown in Fig. \ref{fig:ComparisonClipped}, in all the cases the controller violated at least one of the bounds (100 scenarios by 0.2\,$^\circ$C, 250 scenarios by 0.1\,$^\circ$C), with the controller based on five scenarios violating both the lower and upper bound (by 1.1\,$^\circ$C). Even though the violation decreased with increasing the number of scenarios, further increasing the number of random scenarios to 600 proved unsuccessful in avoiding the violation. Larger problems could not be solved with the given computer in a reasonable amount of time.

If we were to take only extreme values for every uncertainty and consider all the scenarios, we would need to solve a problem with $2^{14+3\times192}$ scenarios, which is too many. To reduce the number of scenarios, we chose to use the nominal scenario, combined with two extreme scenarios. The extreme scenarios were taken as all uncertainties on their lower or upper bound simultaneously. The results of validating the controller for 500 scenarios are shown in the plot `Nominal+two extreme' in Fig.~\ref{fig:Comparison}. The controller based on extreme scenarios was also unable to satisfy the constraints. The black lines after 24 hours cross the green lines so that the lower bound on the temperature is violated (by 0.5$^\circ$\,C). 

\section{Conclusions}
\label{sec:Conclusions}

Solving robust nonlinear optimal control problems is challenging, especially if the knowledge about the uncertainty is limited. Scenario-based approaches provide a way of reformulating the optimal control problems as nonlinear optimization problems. However, the choice of scenarios and their number is non-trivial, because the scenarios must ensure robustness while keeping the size of the optimization problem manageable. In particular, the size of the resulting optimization problems increases computationally if time-varying uncertainties are involved. In this paper, we showed that a class of robust optimal control problems is equivalent to semi-infinite optimization problems. We then demonstrated how a local reduction method derived from semi-infinite optimization provides flexibility in choosing scenarios for nonlinear robust control problems. By adding interim worst-case scenarios, the local reduction method enables finding a trade-off between the size of the resulting optimization problem and robustness of the solution to the original optimal control problem. 

This paper extended the original local reduction method to optimal control problems with time-varying uncertainty. The performance of our approach was evaluated in a case study with both additive and parametric uncertainty. A comparison with common approaches based on random choice of scenarios and on boundary scenarios indicates that local reduction has potential for solving robust optimal control problems in an efficient way while ensuring robustness. 

Future work will include theoretical analysis of the local reduction method as well as numerical improvements needed if local solvers are used. Furthermore, applications in predictive control settings, including analysis of time to find a solution, can also be considered.

\balance
\bibliographystyle{abbrv}
\bibliography{bibTAC}             
                                                   
\end{document}